\documentclass[letterpaper, 10 pt, conference]{ieeeconf} 
\IEEEoverridecommandlockouts 
\usepackage{amsmath,amsfonts,amssymb}

\usepackage{amsthm}
\usepackage{comment}
\usepackage{bbm}
\usepackage{hyperref}
\usepackage{orcidlink}
\usepackage{algorithm}
\usepackage{algorithmic} 
\usepackage{booktabs}

\newtheorem{theorem}{Theorem}[section]

\newtheorem{proposition}[theorem]{Proposition}

\theoremstyle{definition}

\newtheorem{remark}[theorem]{Remark}

\newcommand{\norm}[2][]{\left\lVert#2\right\rVert_{#1}}

\newcommand{\CU}{\widetilde{\mathcal{U}}}

\title{\LARGE \bf
Sparse stabilization of mean-field agent dynamics through a three-operator splitting method*
}

\author{Giacomo Albi$^{1}$, Dante Kalise$^{2}$, Chiara Segala$^{3}$ and Franco Zivcovich$^{4}$
\thanks{*This work has been partially supported by MUR project PRIN 2022 PNRR nr.~P2022JC95T, CUP: B53D23027840001; by Swiss National Science Foundation (SNSF) through the grant number 215528; GA and CS are memeber of INdAM GNCS group.}
\thanks{$^{1}$Giacomo Albi is with the Department of Computer Science,
        University of Verona, IT, %
        {\tt\small giacomo.albi@univr.it}}%
\thanks{$^{2}$Department of Mathematics, Imperial College London, UK,
        {\tt\small dkaliseb@imperial.ac.uk}}%
\thanks{$^{3}$Faculty of Informatics, Università della Svizzera italiana, Lugano, CH,
        {\tt\small chiara.segala@usi.ch}}%
\thanks{$^{4}$Franco Zivcovich is with Department of Computer Science,
        University of Verona, IT, %
        {\tt\small franco.zivcovich@univr.it}}
        }

\begin{document}

\maketitle
\thispagestyle{empty}
\pagestyle{empty}

\begin{abstract}
We study the sparse stabilization of nonlinear multi-agent systems within a mean-field optimal control framework. The goal is to drive large populations of interacting agents toward consensus with minimal control effort. In the mean-field limit, the dynamics are described by a Vlasov-type kinetic equation, and sparsity is enforced through an $\ell_1$–$\ell_2$ penalization in the cost functional. The resulting non-smooth optimization problem is solved via a three-operator splitting (TOS) method that separately handles smooth, non-smooth, and constraint components through gradient, shrinkage, and projection steps. A particle-based Monte Carlo discretization with random batch interactions enables scalable computation while preserving the mean-field structure. Numerical experiments on the Cucker–Smale model demonstrate effective consensus formation with sparse, localized control actions, confirming the efficiency and robustness of the proposed approach.

\end{abstract}
\section{INTRODUCTION}
This paper investigates the sparse stabilization of nonlinear, second-order multi-agent systems (MAS), where each agent is characterized by position and velocity and interacts through nonlocal communication laws encoding attraction, repulsion, and alignment. Such interactions generate emergent collective behaviors typical of self-organizing biological systems~\cite{cucker2007IEEE,d2006self,couzin2005N} and appear in diverse contexts, from crowd motion and social dynamics~\cite{hegselmann2002opinion,cristiani2014book,MR3542027,AlbiFerrareseSegala2021} to robotic swarms and AI networks~\cite{zhao2018self,duarte2016evolution,ha2022collective}. 

The collective dynamics of large systems can be described at multiple scales: microscopically by ODEs and at mesoscopic or macroscopic levels by mean-field PDEs~\cite{carrillo2010particle,PT:13,maffettone2023continuification}. Deriving higher-scale models such as mean-field or kinetic equations \cite{AHP15} is crucial both theoretically and computationally. Mean-field formulations consistently capture macroscopic behavior by averaging microscopic fluctuations, while reducing the computational cost of simulating high-dimensional $N$-agent systems~\cite{albi2013binary,jin2020random,albi2025control}.

Our main goal is to design control strategies that steer the system toward consensus, namely a collective state where all agents share a common velocity. As a prototype, we consider the controlled Cucker--Smale model~\cite{CFPT13,BAILO20181}, a benchmark for self-organization and control in collective dynamics. We seek parsimonious controls acting on few agents to minimize intervention cost. Sparsity is promoted by an $\ell_1$-penalty in the cost functional~\cite{candes2006stable,donoho2006compressed,cheng2015handbook,cohen1983absolute,BonginiFornasier2017}, which effectively targets the most misaligned agents~\cite{CFPT15}. Since the non-smooth $\ell_1$ term prevents standard gradient methods, we employ a three-operator splitting (TOS) algorithm~\cite{beck2017first,DavisYin17} that alternates gradient and proximal steps over smooth, non-smooth, and constrained components. Gradients are computed via adjoint equations, while proximal maps admit closed-form shrinkage and projection operators.

To ensure scalability, we adopt a mean-field formulation leading to a Vlasov-type kinetic equation for the agent density in phase space~\cite{FPR14,MR3264236,Lacker2017}. This model preserves the essential microscopic structure while lowering computational cost and is discretized through a Monte Carlo particle scheme with linear complexity~\cite{albi2017mean,PT:13}. Combining mean-field modeling with proximal splitting yields a unified, scalable framework for sparse stabilization of collective dynamics. Numerical experiments on the alignment model demonstrate consensus with minimal control effort and persistent sparsity in both time and across agents.

The paper is organized as follows:
Section~\ref{sec:mf} introduces the mean-field  optimal control formulation from microscopic to mesoscopic level;
Section~\ref{sec:particle_based_scheme}
details the particle discretiazion of the mean-field control problem; Section~\ref{sec:prox} describe the TOS gradient for the mean-field sparse stabilitzation. Finally,
Section~\ref{sec:numerics} presents different numerical results.

\section{MEAN-FIELD APPROXIMATION OF OPTIMAL MAS CONTROL}\label{sec:mf}
We consider a MAS of $N$ interacting agents with states $(x_i(t), v_i(t)) \in \mathbb{R}^d \times \mathbb{R}^d$, where $d$ denotes the dimension of the physical space. The evolution of each agent is governed by the second-order dynamics
\begin{equation}\label{eq:CS}
	\begin{split}
	\dot{x}_i &= v_i ,\\
	\dot{v}_i &= \frac{1}{N} \sum_{j=1}^N P(\|x_i-x_j\|_2)(v_j - v_i) + u_i , 
	\end{split}
\end{equation}
for $i = 1, \ldots, N$, with initial data $x_i(0) = x_i^0$ and $v_i(0) = v_i^0$. The function $P(\cdot)$ is a symmetric, radial communication kernel describing pairwise interactions between agents, defined as
\begin{align}\label{eq:Prad}
P(r) =(1+r^2)^{-\kappa},\qquad \kappa>0.
\end{align}
Here, the objective is to design control laws $u_i(t)$ that drive the system toward a state of consensus where all agents share the same velocity
$\bar{v}(t) = (N)^{-1}\sum_{j=1}^N v_j(t).$ The term $u(t) = (u_1(t), \ldots, u_N(t))\in\mathbb{R}^{dN}$ denotes the collection of control inputs acting on the agents, where each $u_i: [0,T] \to U \subset \mathbb{R}^{d}$ and $U$ a compact set of admissible controls, moreover, in what follows we will require that the set of admissible control is such that 
\begin{equation}\label{eq:Uadm}\mathcal{U}^N=\left\{u(\cdot): \frac 1 N\int_0^T\sum_{i=1}^N\|u_i(\tau)\|_1 \,d\tau\leq B\right\},\end{equation}
where $B>0$ is a fixed budget. We assume that the control vector $u(\cdot)$ is obtained as 
 minimizer over $\mathcal{U}^N$ of a cost functional $\mathcal{J}^N(u)\equiv \mathcal{J}^N(u;x^0,v^0)$ defined by
\begin{equation}
\label{eq:minCS}
\begin{split}
	\mathcal{J}^N(u)=
	&\int_0^T \frac{1}{N} \sum_{i=1}^N
	\Big(
	\| \bar{v}(t) - v_i(t) \|_2^2
	\cr
     &\quad+ \alpha \|u_i(t)\|_2^2
	+ \beta \|u_i(t)\|_1
	\Big) \mathrm{d}t,
    \end{split}
\end{equation}
where $\alpha, \beta \ge 0$ are penalization parameters. The first term promotes consensus among the agents, the quadratic term penalizes the overall control energy, and the $\ell_1$-term introduces sparsity, ensuring that only a limited number of agents or time instants are actively controlled. This type of formulation has been widely employed in the context of sparse control of alignment dynamics \cite{CFPT13,bongini2014sparse,kalise2017infinite}.

Although the microscopic control problem \eqref{eq:CS}–\eqref{eq:minCS} provides an accurate description of the dynamics, its computational cost grows rapidly with $N$. To overcome this limitation, a mean-field approximation can be derived by considering the limit $N \to \infty$. In this limit, the individual-based model is replaced by a continuous distribution of agents described by a probability density function $f(t,x,v)$ defined over the phase space $\mathbb{R}^{2d}$. The corresponding empirical measure associated with the microscopic system is given by
$f^N(t,x,v) = {N}^{-1}\sum_{i=1}^N \delta(x - x_i(t))\otimes\,\delta(v - v_i(t)).$
As $N \to \infty$, the empirical measure $f^N$ converges to $f(t,x,v) \in \mathcal{M}_1(\mathbb{R}^{2d})$, the space of probability measures with compact support, see \cite{MR4028474,MR3264236}. The evolution of $f$ is governed by a Vlasov-type  equation of the form
\begin{equation}\label{eq:mean-field}
	\partial_t f + v \cdot \nabla_x f = - \nabla_v \cdot \big[ f\,(\mathcal{P}[f] + u) \big],
\end{equation}
with initial data $ f(0,x,v) = f^0(x,v),$ and  where the nonlocal alignment operator $\mathcal{P}[f]\equiv \mathcal{P}[f](t,x,v)$ is
\[
\mathcal{P}[f](t,x,v) = \int_{\mathbb{R}^{2d}} P(\|x - x_*\|_2)(v_* - v)\, df_*,
\]
where we used the shorten notation $df_*=f(t,x_*,v_*)\textrm{d}x\textrm{d}v$.
The control $u = u(t,x,v)$ is a distributed control acting such that
$ u:[0,T]\times \mathbb{R}^{2d} \to \mathbb{R}^d,$ and obtained by minimizing the following cost functional
\begin{equation}\label{eq:cost_mf}
	\begin{split}
    \mathcal{J}(u) :=
	&\int_0^T \int_{\mathbb{R}^{2d}}
	\Big(
	\|m[f](t) - v\|_2^2
    \cr
	&
    + \alpha \|u\|_2^2
	+ \beta \|u\|_1
	\Big)	df\,\mathrm{d}t,
    \end{split}
\end{equation}
where
\[
m[f](t) = \int_{\mathbb{R}^{2d}} v\,df(x,v),
\]
is the mean velocity of the population, and the set of the admissible control in the mean-field limit is denoted as $\mathcal{U}$, where the budget constrain in \eqref{eq:Uadm} is replaced by
\begin{equation}\label{eq:Uadm_mf}
\int_0^T\int_{\mathbb{R}^{2d}}\|u(\tau,x,v)\|_1 \,d f\,d\tau\leq B.
\end{equation}

\section{PARTICLE-BASED SCHEME FOR MEAN-FIELD SPARSE STABILIZATION}\label{sec:particle_based_scheme}
To efficiently approximate the mean-field control dynamics \eqref{eq:mean-field}--\eqref{eq:cost_mf}, we develop a particle-based methods based on randomized sampling similarly to methods proposed in \cite{albi2013binary,PT:13,jin2020random}. Here, we combine these strategies with optimization methods for the synthesize of the optimal control.

We account for $N_s$ particles $\{(x_i^0, v_i^0)\}_{i=1}^{N_s}$ sampled independently from the initial distribution $f^0(x,v)$. The mean-field dynamics is then approximated by a time-discrete particle system evolving according to the microscopic dynamics \eqref{eq:CS}. We introduce the discrete time grid $t^n = n \Delta t$, with $N_T \Delta t = T$, and simulate the system for $i=1,\ldots,N_s$ and $n=0,\ldots,N_T-1$ as
\begin{equation}\label{eq:disc_particles}
	\begin{split}
		&{x}^{n+1}_i = x^{n}_i + \Delta t \, v_i^n \, , \\
		&{v}^{n+1}_i = v^{n}_i  + \Delta t \, u^n_i \cr
        &\quad+ \frac{\Delta t}{M_s} \sum_{k=1}^{M_s} {P}\!\left(\norm[2]{x_i^n-x_{j_k}^n}\right) \big(v^n_{j_k}-v^n_i\big), \, 
	\end{split}
\end{equation}
where $\{j_1,\ldots,j_{M_s}\}\subseteq\{1,\ldots,N_s\}$ is a random subsample of size $M_s \leq N_s$ drawn at each time step to approximate the interaction operator. Notice that the interaction is computed with respect to a sampled subset of agents, this reduces the cost of each update from $\mathcal{O}(N_s^2)$ to $\mathcal{O}(N_s M_s)$, making the scheme feasible for large-scale simulations.

The performance of the control strategy is assessed by discretizing the mean-field cost functional \eqref{eq:cost_mf}. The discrete counterpart reads as follows
\begin{equation}\label{eq:J_discrete}
\begin{split}
	&\mathcal{J}_{\Delta t}^{N_s}({\boldsymbol{u}}) :=
	\frac{\Delta t}{N_s}
    \sum_{n=0}^{N_T-1}\sum_{j=1}^{N_s} 
	\left(
	\norm[2]{ \bar{v} - v^n_j }^2\right.\cr
	&\qquad\left.+ \alpha \norm[2]{ u_j^n}^2 
	+ \beta \norm[1]{ u_j^n}
	\right) , 
    \end{split}
\end{equation}
where $\{u^n\}=(u_1^n,\ldots,u_{N_s}^n)\in\mathbb{R}^{dN_s}$ acts on each interval $[t^n,t^{n+1}]$, and 
we introduce the concatenated control vector
\[
\boldsymbol{u} := \big( u^0, \, u^1 \, \ldots\, u^{N_T-1} \big) \in \mathbb{R}^{dN_sN_T},
\]
which collects all the agents' controls over the discretized time horizon. The admissible set of control \eqref{eq:Uadm_mf} in the discrete setting corresponds to 
\begin{align}\label{eq:Uadm_disc}
\CU = \Big\{ \boldsymbol{u} \in \mathbb{R}^{d N_s N_T} : \norm[1]{\boldsymbol{u}} \le \widetilde{B} \Big\},
\end{align}
which is a non-empty, closed, and convex set, and where $\widetilde B=N_s B/\Delta t>0$ denotes the scaled budget according to \eqref{eq:Uadm_mf}. The discrete optimization problem reads

 \begin{equation}
 \begin{aligned}\label{eq:opt_disc}
 &\min_{\boldsymbol{u}\in \CU } \mathcal{J}_{\Delta t}^{N_s}(\boldsymbol{u}),\quad \textrm{s.t.}\, \eqref{eq:disc_particles}.
 \end{aligned}
 \end{equation}

\subsection{Discrete optimality conditions} 
We derive  from \eqref{eq:disc_particles}-\eqref{eq:J_discrete} the set of optimality conditions for the discrete agent system introducing a Lagrangian functional that enforces the discrete dynamics \eqref{eq:disc_particles} as constraints. Hence,
let ${ p}^n=(p_1^n,\ldots,p_{N_s}^n)$ and ${ q}^n=(q_1^n,\ldots,q_{N_s}^n)$ be the adjoint multipliers associated with the state variables ${ x}^n$ and ${ v}^n$, respectively, then 
the adjoint dynamics reads
\begin{equation}\label{eq:disc_adjoint}
\begin{aligned}
	&p_i^n = p_i^{n+1} +  \frac{ \Delta t}{M_s}\sum_{k=1}^{M_s} P'\!\left(\|x_i^n-x_{j_k}^n\|_2\right)\cr
    &\quad\langle q^{n+1}_{j _k}- q^{n+1}_{i}, v^n_{j_k} - v^n_{i}\rangle (x^n_{j_k}- x^n_{i}) , \\
	&q_i^n = q_i^{n+1} + 2\Delta t\big(v_i^n - \bar{v}^n\big)- \Delta tp_i^{n+1}\cr
    &\,+ \frac{\Delta t}{M_s}\sum_{k=1}^{M_s} P\!\left(\|x_i^n-x_{j_k}^n\|_2\right)(q_i^{n+1}-q_{j_k}^{n+1}) ,
\end{aligned}
\end{equation}
with terminal conditions $p_i^{N_T}=0$, $q_i^{N_T}=0$. 
In the unconstrained setting, the variation with respect to $\boldsymbol{u}$ yields, for each $i=1,\dots,N_s$ and $n=0,\dots,N_T-1$, 
to the optimality condition
\begin{equation}\label{eq:oc_discrete_time}
	0 \in \bigl( 2\alpha u_i^n + q_i^{n+1} + \beta\,\partial\|u_i^n\|_1 \bigr),
\end{equation}
where $\partial\|\cdot\|_1$ denotes the subdifferential of the $\ell_1$-norm. Under the admissibility constraints $u_i^n \in \mathcal{U}_{\mathrm{adm}}$, 
the condition is complemented by the normal cone to the admissible set.
\begin{remark}\label{rmk:micro_mf}
We remark that when $N_s = M_s = N$ the proposed particle discretization coincides to the time-discrete formulation of the optimality system associated to the optimal control problem for the microscopic system \eqref{eq:CS}–\eqref{eq:minCS}, and similarly to \cite{BAILO20181,CFPT13}. For further details about the consistency results of the discretized adjoint system and its mean-field formulation we refer to \cite{albi2017mean, herty2019consistent,burger2021mean}.
\end{remark}
\section{TOS scheme for mean-field sparse control}\label{sec:prox}
We address the solution of the discrete optimization problem \eqref{eq:opt_disc}, which is a high-dimensional non-smooth optimization problem due to the $\ell_1$-term promoting sparsity in the control vector $\boldsymbol{u}\in\mathbb{R}^{dN_sN_T}$ and to the presence of a non-linear interaction kernel. To efficiently handle this structure, we adopt a three-operator splitting (TOS) algorithm, which extends proximal gradient methods to combine smooth quadratic terms, the non-smooth $\ell_1$ penalty, and possible convex constraints. The solution of the resulting optimality system \eqref{eq:disc_particles}–\eqref{eq:disc_adjoint}–\eqref{eq:oc_discrete_time} is then obtained through the iterative TOS procedure.
The discrete optimization problem \eqref{eq:opt_disc} can be rewritten in the equivalent TOS form 
\begin{equation}\label{eq:ista_pb}
\min_{\boldsymbol{u}\in\mathbb{R}^{dN_sN_T}}\mathcal{J}_1(\boldsymbol{u}) + \mathcal{J}_2(\boldsymbol{u}) + \mathcal{J}_3(\boldsymbol{u}) ,
\end{equation}
where the first two components correspond respectively to the smooth term
\begin{align}\label{eq:j1_disc}
\mathcal{J}_1(\boldsymbol{u}) =  \frac{\Delta t}{ N_s}\sum_{n=0}^{N_T-1} \sum_{i=1}^{N_s} 
\Big( \norm[2]{\bar{v}^n - v_i^n}^2 + \alpha \norm[2]{u_i^n}^2 \Big),
\end{align}
and the non-smooth parts related to the $\ell_1$ penalty and the budget constraint as follows  
\[
\mathcal{J}_2(\boldsymbol{u}) =  \frac{\Delta t}{ N_s}\sum_{n=0}^{N_T-1} \sum_{i=1}^{N_s} \beta \norm[1]{u_i^n}  =  \frac{\Delta t}{N_s}\beta \norm[1]{\boldsymbol{u}},
\]
 and the barrier function
\[
\mathcal{J}_3(\boldsymbol{u}) = \delta_{\CU}(\boldsymbol{u}) :=
\begin{cases}
	0, & \text{if } \boldsymbol{u} \in \CU, \\
	+\infty, & \text{if } \boldsymbol{u} \notin \CU,
\end{cases}
\]
where the admissible set is defined in \eqref{eq:Uadm_disc}.
Following the approach of \cite{DavisYin17}, we employ the TOS scheme as a generalization of the classical gradient method, incorporating two proximal mappings for the non-smooth part, and corresponding to the shrinkage and projection operators. Starting from an initial guess $\boldsymbol{u}^{(0)} \in \mathbb{R}^{d N_s N_T}$, the algorithm generates an iterative sequence ${\boldsymbol{u}^{(k)}}$ via
\begin{equation} \label{eq:tos}
	\begin{split}
	\boldsymbol u^{(k)}_{\mathcal{J}_3} &= \text{prox}_{\lambda \mathcal{J}_3} \left(\boldsymbol{u}^{(k)} \right), \\
	\boldsymbol u^{(k)}_{\mathcal{J}_2} &= \text{prox}_{\lambda \mathcal{J}_2} \left( 2 \boldsymbol u^{(k)}_{\mathcal{J}_3} - \boldsymbol{u}^{(k)}\right.\cr
    &\qquad\quad\qquad\left.- \lambda  \nabla \mathcal{J}_1 \left(\boldsymbol u^{(k)}_{\mathcal{J}_3}\right)  \right), \\
	\boldsymbol{u}^{(k+1)} &= \boldsymbol{u}^{(k)}  + \lambda_k \left( \boldsymbol u^{(k)}_{\mathcal{J}_2} -  \boldsymbol u^{(k)}_{\mathcal{J}_3} \right).
		\end{split}
\end{equation}
Since $\mathcal{J}_1$ depends on both $\boldsymbol{u}$ and the velocity variables $\boldsymbol{v}$ through the discrete dynamics \eqref{eq:disc_particles}, its gradient with respect to $\boldsymbol{u}$ must account for this coupling. Using the adjoint variables $q_i^{n+1}$ to represent this dependence, the gradient is
\begin{equation}\label{eq:opt_J1}
\nabla \mathcal{J}_1(u_i^n) = \frac{\Delta t}{ N_s}(2\alpha u_i^n + q_i^{n+1}),
\end{equation}
which gives the discrete-time optimality condition as the smooth component of \eqref{eq:oc_discrete_time}. The parameter $\lambda \in (0, 2\vartheta)$ now plays the role of the stepsize for the forward (gradient) step, consistent with the notation in \cite{DavisYin17}, where $\vartheta > 0$ denotes the cocoercivity constant of $\nabla \mathcal{J}_1$, inversely related to its Lipschitz constant.
The relaxation parameter $\lambda_k \in \bigl(0, (4\vartheta - \lambda)/(2\vartheta)\bigr)$ controls the averaging between successive iterates.

We introduced in \eqref{eq:tos} the proximal operator, i.e. for a given convex, lower semi-continuous function $\varphi$ defined as
\[
\text{prox}_{\lambda \varphi}(w) = \arg\min_{z} \left\{ \frac{1}{2} \| z - w \|_2^2 + \lambda \varphi(z) \right\}.
\]
To explicitly compute the first and second step of \eqref{eq:tos} we can rely on the following proposition

\begin{proposition}\label{prop:prox_ops}
The proximal operators of $\lambda \mathcal{J}_2$ and $\lambda \mathcal{J}_3$ in \eqref{eq:tos} are given respectively by
\begin{align}
\text{prox}_{\lambda \mathcal{J}_2}(w)
&= \mathbb{S}_{\beta\lambda \frac{\Delta t}{ N_s}}(w)\label{eq:softr}\\
\text{prox}_{\lambda \mathcal{J}_3}(w)
&= \mathbb{P}_{\CU}(w),\label{eq:proj}
\end{align}
where the soft-thresholding operator is defined as
\[
\mathbb{S}_{\textrm{h}}(w):= \operatorname{sign}(w)\,\left[\,|w| - \textrm{h}\,\right]_+,
\]
with $[\,\cdot\,]_+$ the positive part applied componentwise,
and the projection map is
\[
\mathbb{P}_{\CU}(w) :=
\begin{cases}
w, & \text{if } \|w\|_1 \le \widetilde B, \\[4pt]
\mathbb{S}_{\lambda^*}(w), & \text{if } \|w\|_1 > \widetilde B,
\end{cases}
\]
with $\lambda^*>0$ chosen such that $\|\mathbb{S}_{\lambda^*}(w)\|_1 = \widetilde B$.
\end{proposition}
These results follow from standard proximal operator derivations; see \cite[Ch.~6, Sec.~6.2--6.4]{beck2017first} for details.

\begin{remark}

The successive proximal steps in \eqref{eq:tos} ensure stability and convergence.  
The soft-thresholding operator $\mathbb{S}(\cdot)$ performs the ISTA-type shrinkage that promotes sparsity, while the projection $\mathbb{P}(\cdot)$ enforces feasibility by constraining the control to $\CU$.
\end{remark}
\begin{algorithm}[!t]
\caption{TOS reduced gradient}
\label{alg:TOS}
\begin{algorithmic}[1]
  \STATE \textbf{Input:} 
  Inital density $f^0(x,v)$ and final time $T>0$, time step $\Delta t>0$, number of time steps $N_T$, number of sampled agents $N_s \ge M_s$, 
  parameters $\alpha,\beta>0$ (see~\eqref{eq:J_discrete}), 
  tolerance $tol>0$, maximum iterations $k_{\max}$, 
  and initial guess $\boldsymbol{u}^{(0)}$.
  \STATE Sample initial states $\{(x_i^0,v_i^0)\}_{i=1}^{N_s}$ from $f^0(x,v)$.
  \STATE Set iteration counter $k=0$.
  \WHILE{$\|\nabla \mathcal{J}(\boldsymbol{u}^{(k)})\|_2 > tol$ \textbf{and} $k<k_{\max}$}
        \STATE (i)~Compute $(\boldsymbol{x}^{(k)},\boldsymbol{v}^{(k)})$ via \eqref{eq:disc_particles}.
        \STATE (ii)~Compute $(\boldsymbol{p}^{(k)},\boldsymbol{q}^{(k)})$ via \eqref{eq:disc_adjoint}.
    \STATE (iii)~Set step parameters $\lambda,\,\lambda_k>0$.
    \STATE (iv)~Projection step $\boldsymbol{u}^{(k)}_{\mathcal{J}_3}$ via \eqref{eq:proj}.
    \STATE (v)~Evaluate gradient $\nabla \mathcal{J}_1(\boldsymbol{u}^{(k)}_{\mathcal{J}_3})$ as in \eqref{eq:opt_J1}.
    \STATE (vi)~Compute auxiliary variable:
    \[
      \boldsymbol{w}^{(k)} = 2\boldsymbol{u}^{(k)}_{\mathcal{J}_3} - \boldsymbol{u}^{(k)} - \lambda \nabla \mathcal{J}_1(\boldsymbol{u}^{(k)}_{\mathcal{J}_3}).
    \]
    \STATE (vii)~Shrinkage step $\boldsymbol{u}^{(k)}_{\mathcal{J}_2} = \text{prox}_{\lambda \mathcal{J}_2}(\boldsymbol{w}^{(k)})$ via \eqref{eq:softr}.
    \STATE (viii)~Update control:
    \[
      \boldsymbol{u}^{(k+1)} = \boldsymbol{u}^{(k)} + \lambda_k \big(\boldsymbol{u}^{(k)}_{\mathcal{J}_2} - \boldsymbol{u}^{(k)}_{\mathcal{J}_3}\big).
    \]
    \STATE (ix)~Increment $k \gets k+1$.
  \ENDWHILE
\end{algorithmic}
\end{algorithm}

\section{Numerical experiments}\label{sec:numerics}
We test the performance of the sparse control strategies at both microscopic and mean-field scales.

\subsection{Test 1: Microscopic dynamics}
We first assess the methodology on the microscopic dynamics \eqref{eq:CS}, as justified by Remark~\ref{rmk:micro_mf}. 
We consider the one-dimensional system of $N=20$ agents with Cucker-Smale dynamics \eqref{eq:CS} with $\kappa =1$ in \eqref{eq:Prad}. The final horizon is set $T=15$ with time discretizion step $\Delta t=0.2$. 
Initial states are sampled as $x_i^0 \sim \mathcal{N}(\mu_x^{(i)}, \sigma_x^2)$, $v_i^0 \sim \mathcal{N}(\mu_v^{(i)}, \sigma_v^2)$,
where $(\mu_x^{(1)},\mu_v^{(1)})=(0,-2)$, $(\mu_x^{(2)},\mu_v^{(2)})=(0,0)$, $(\mu_x^{(3)},\mu_v^{(3)})=(-2,2)$, 
with $\sigma_x=0.2$, $\sigma_v=0.4$, forming three clusters that prevent consensus without control~\cite{cucker2007IEEE}.
We consider the control functional \eqref{eq:minCS} with a fixed budget as in \eqref{eq:Uadm_disc}, setting $\widetilde B = 120$ and $\alpha = 0.01$ for the $\ell_2$ regularization term, while varying $\beta$ to induce sparsity through the $\ell_1$ penalty.
The sparse control problem is solved using the TOS scheme~\eqref{eq:tos} with the following fixed parameters: 
$tol = 10^{-8}$, $k_{\max} = 5\times 10^2$, $\lambda= 0.1$, $\lambda_k = 1$, and $\boldsymbol{u}^{(0)}=0$. The 
We test $\beta \in \{5\times10^{-2}, 10^{-1}, 1\}$ to study the effect of sparsity.
\begin{figure}[h!]
	\centering
	\includegraphics[width=0.485\linewidth]{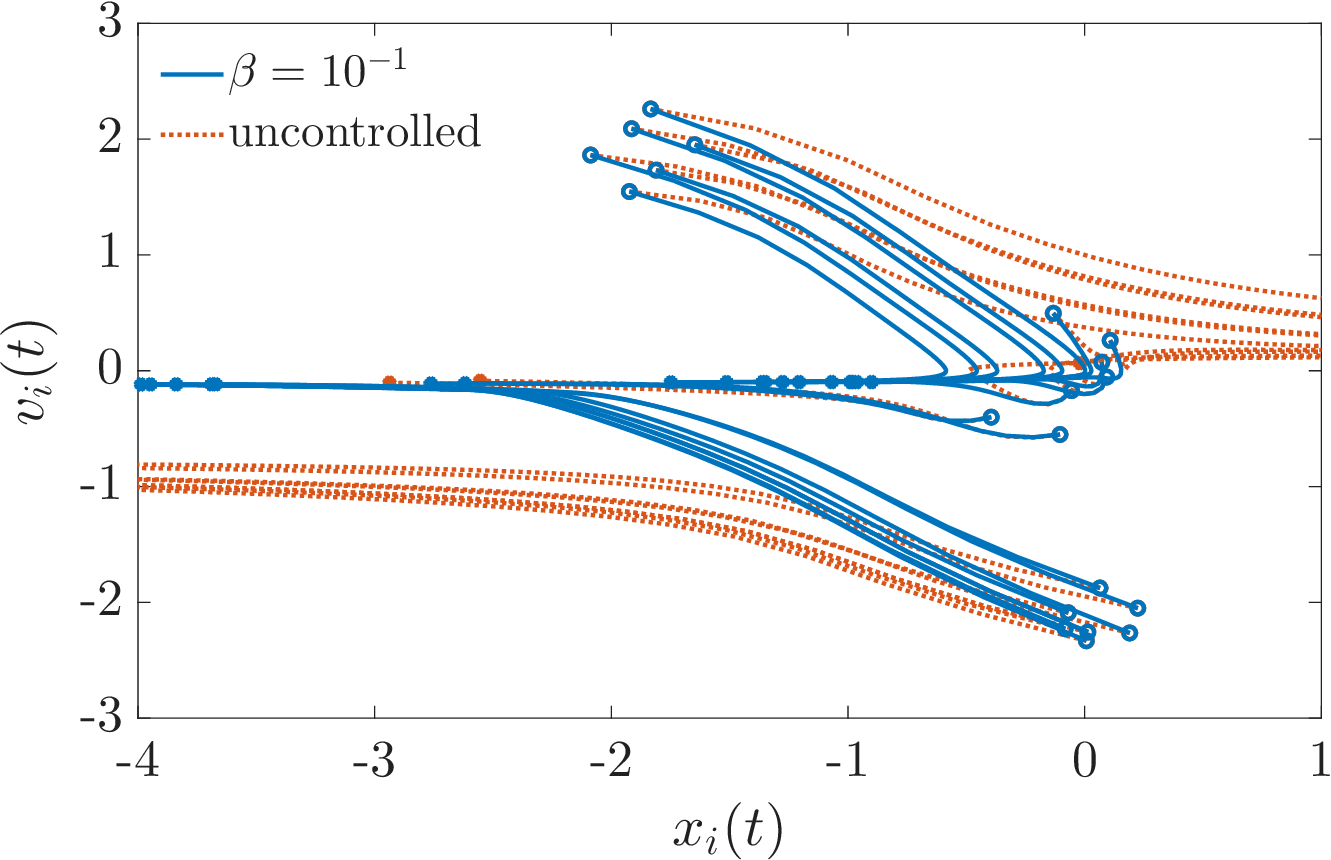}
\includegraphics[width=0.485\linewidth]{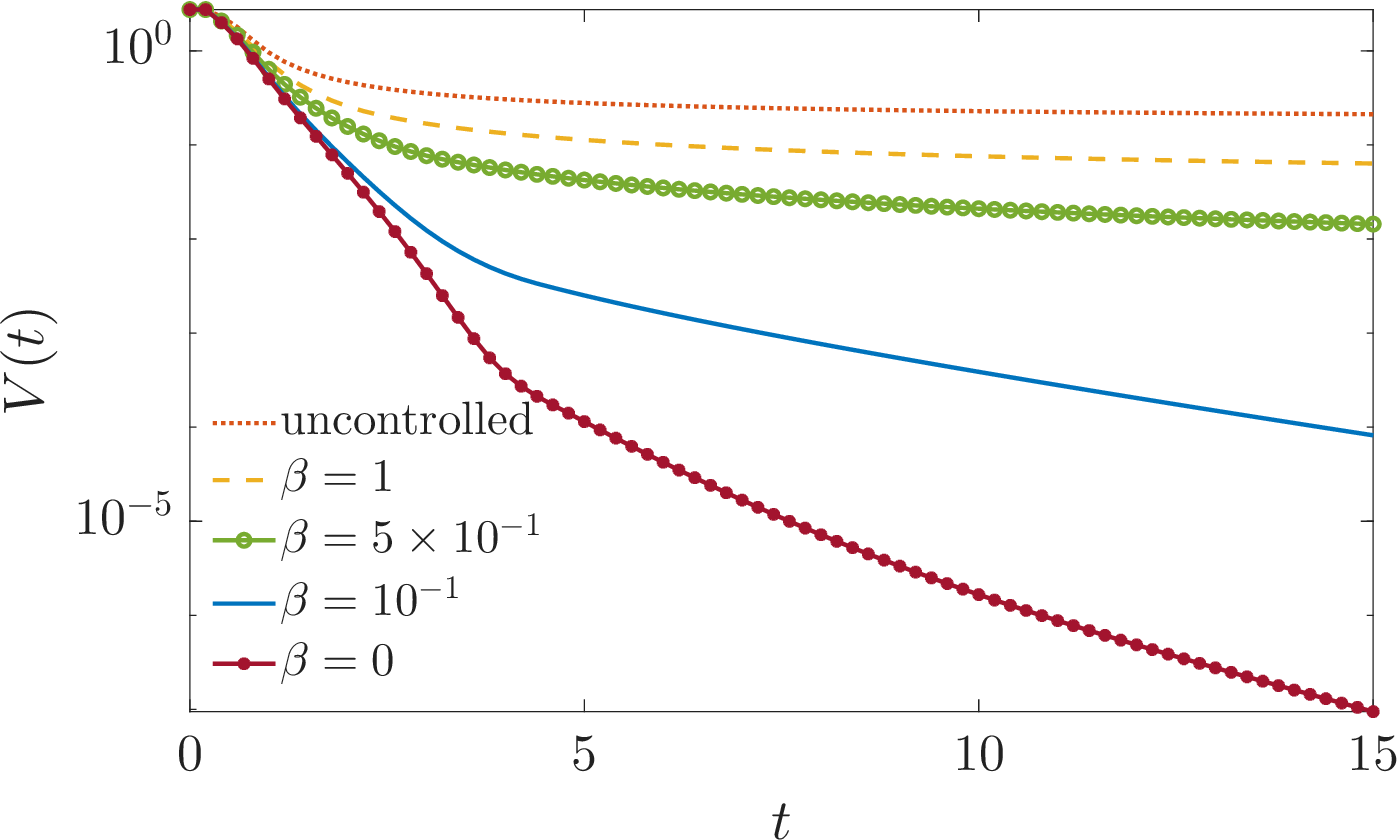}\\[2pt]
\includegraphics[width=0.485\linewidth]{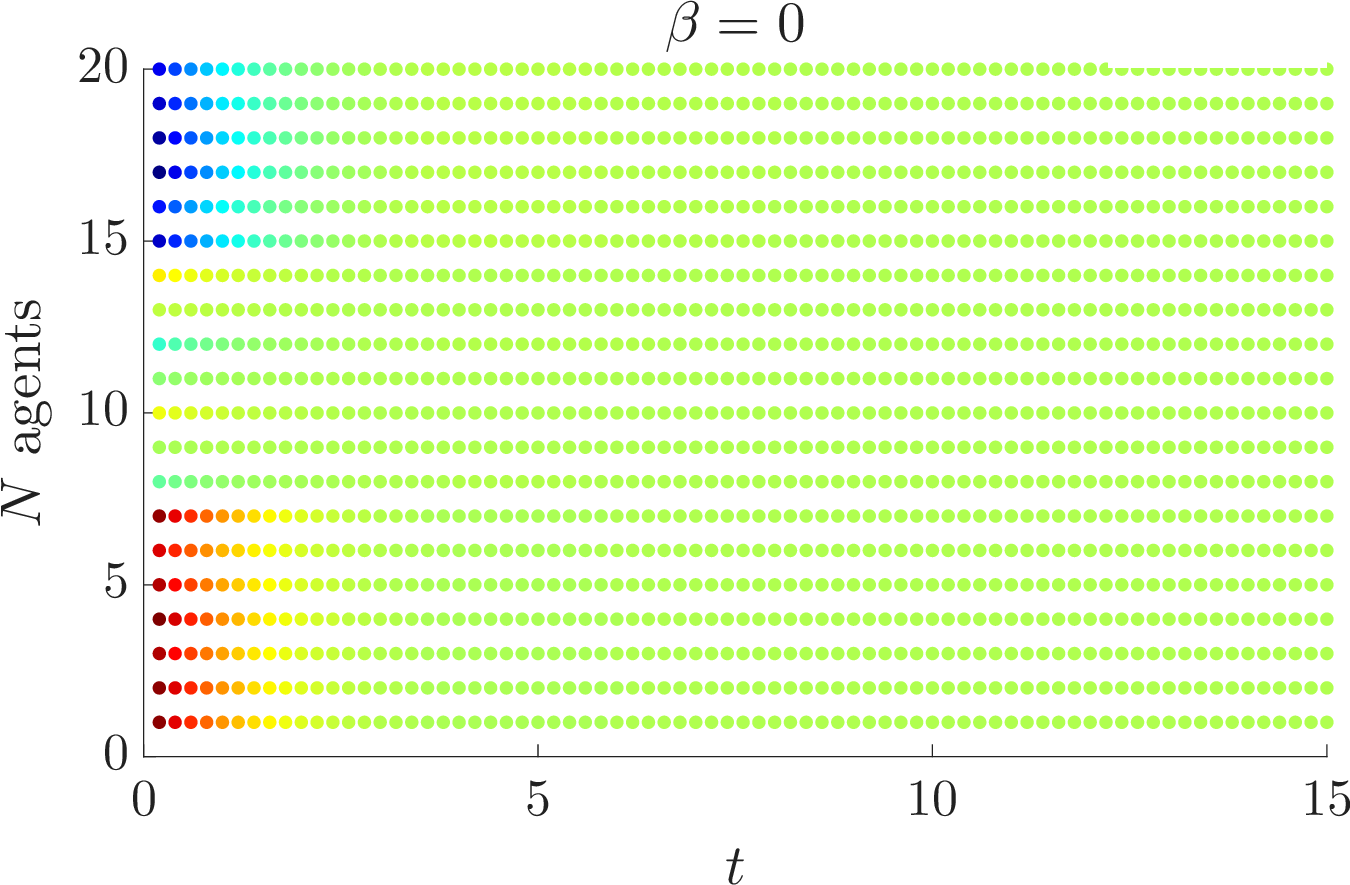}
	\includegraphics[width=0.485\linewidth]{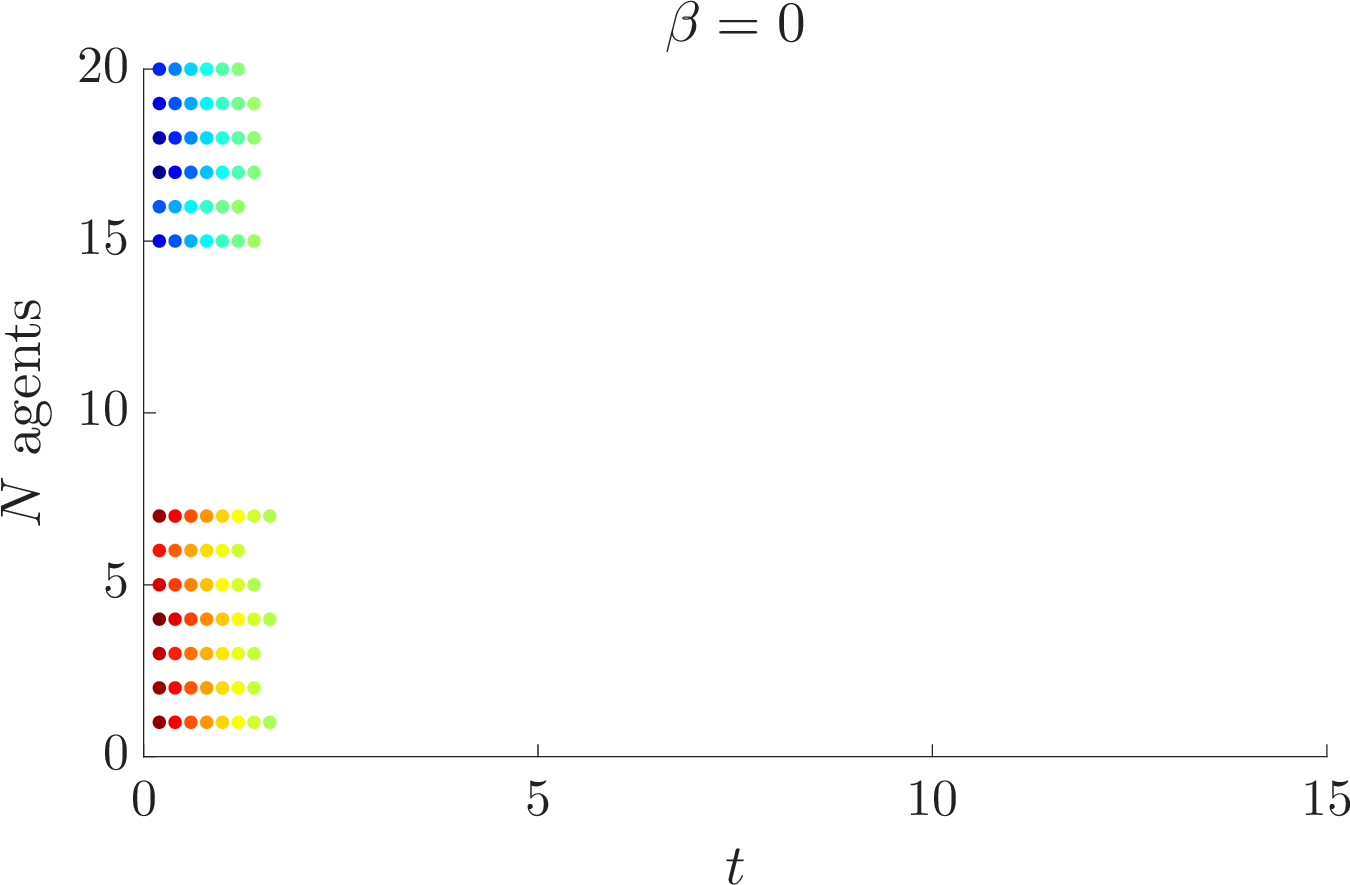}
	\caption{\textbf{Test 1.} \emph{Top row:} (left) uncontrolled vs.\ controlled trajectories for $\alpha=10^{-2}$, $\beta=10^{-1}$; (right) decay of the Lyapunov functional $V(t)={(2N)^{-2}}\sum_{i,j}|v_i(t)-v_j(t)|^2$ for different $\beta$. 
 \emph{Bottom row:} (left) active control components for $\beta = 0$ and (right) for $\beta=10^{-1}$.}
	\label{fig:uncontrolled_vs_controlled}
\end{figure}

Figure~\ref{fig:uncontrolled_vs_controlled} compares uncontrolled, sparse controlled, and optimally controlled dynamics. 
For $\beta=10^{-1}$, agents reach consensus while the uncontrolled system remains dispersed. 
Smaller $\beta$ values accelerate convergence with more active controls, whereas larger $\beta$ increase sparsity but slow alignment, showing a trade-off between efficiency and parsimony.

\subsection{Test 2: Mean-field sparse stabilization}
We next evaluate the sparse control strategy on the mean-field Cucker–Smale model~\eqref{eq:mean-field}. 
The initial density $f^0(x,v)$ reproduces the microscopic configuration of Test~1, consisting of three clusters in phase space centered at $(\mu_x^{(i)},\mu_v^{(i)})$ with spreads $\sigma_x=0.2$ and $\sigma_v=0.4$. 
The system is initialized with $N_s=4\times10^4$ particles and simulated up to $T=15$ using a time step $\Delta t=0.2$. 
The dynamics are approximated through the Monte Carlo particle scheme described in Section~\ref{sec:particle_based_scheme} with $M_s=100$ random batches.  The mean-field cost functional~\eqref{eq:cost_mf} adopts the same regularization weight $\alpha=0.01$ as in Test~1, while varying the sparsity parameter $\beta$ to assess the influence of the $\ell_1$ penalization. 
The TOS algorithm is applied with the same numerical settings and stopping criteria used in the microscopic case.

\begin{figure}[htbp]
    \centering
    \includegraphics[width=0.485\linewidth]{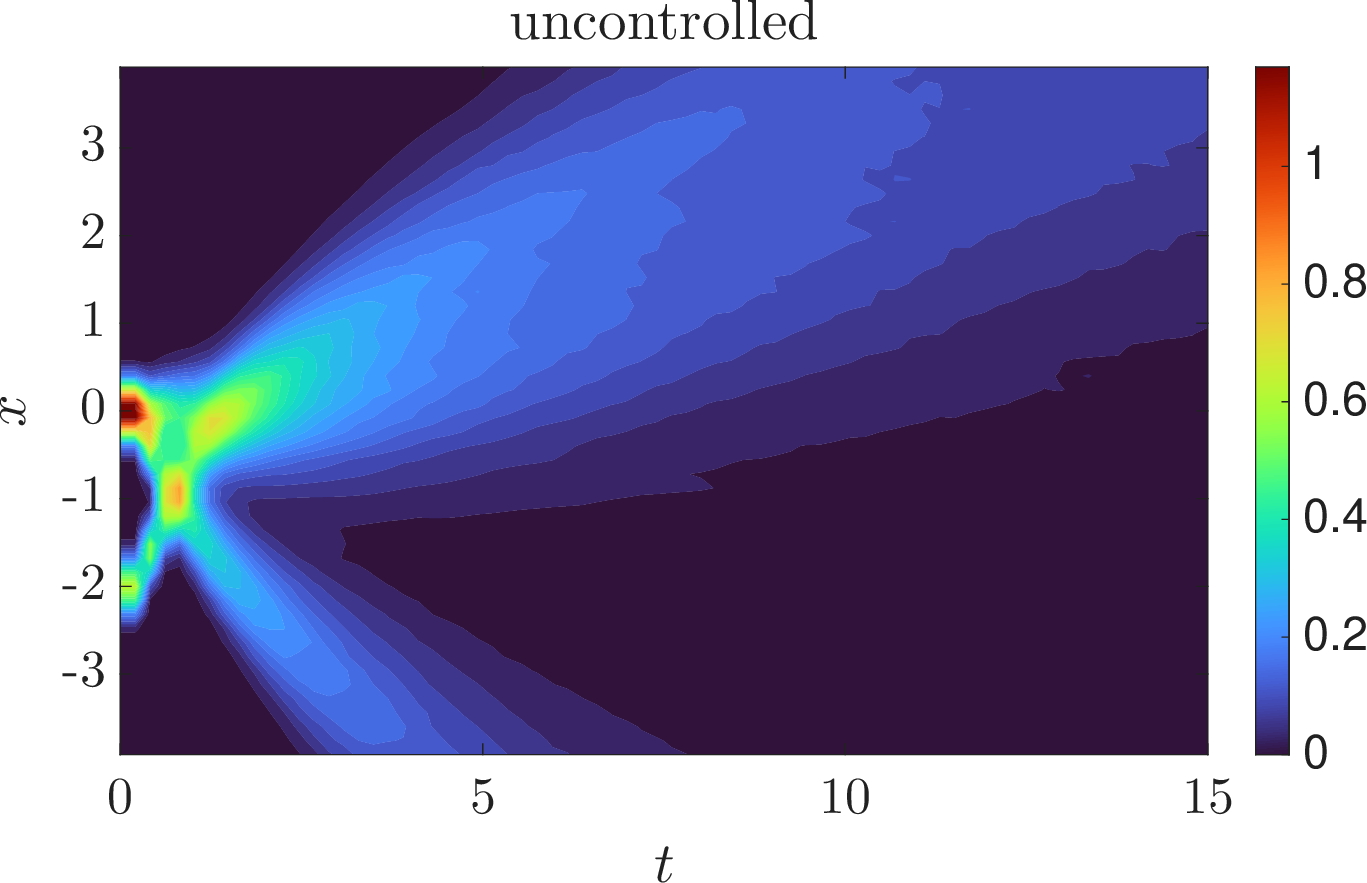}
    \includegraphics[width=0.485\linewidth]{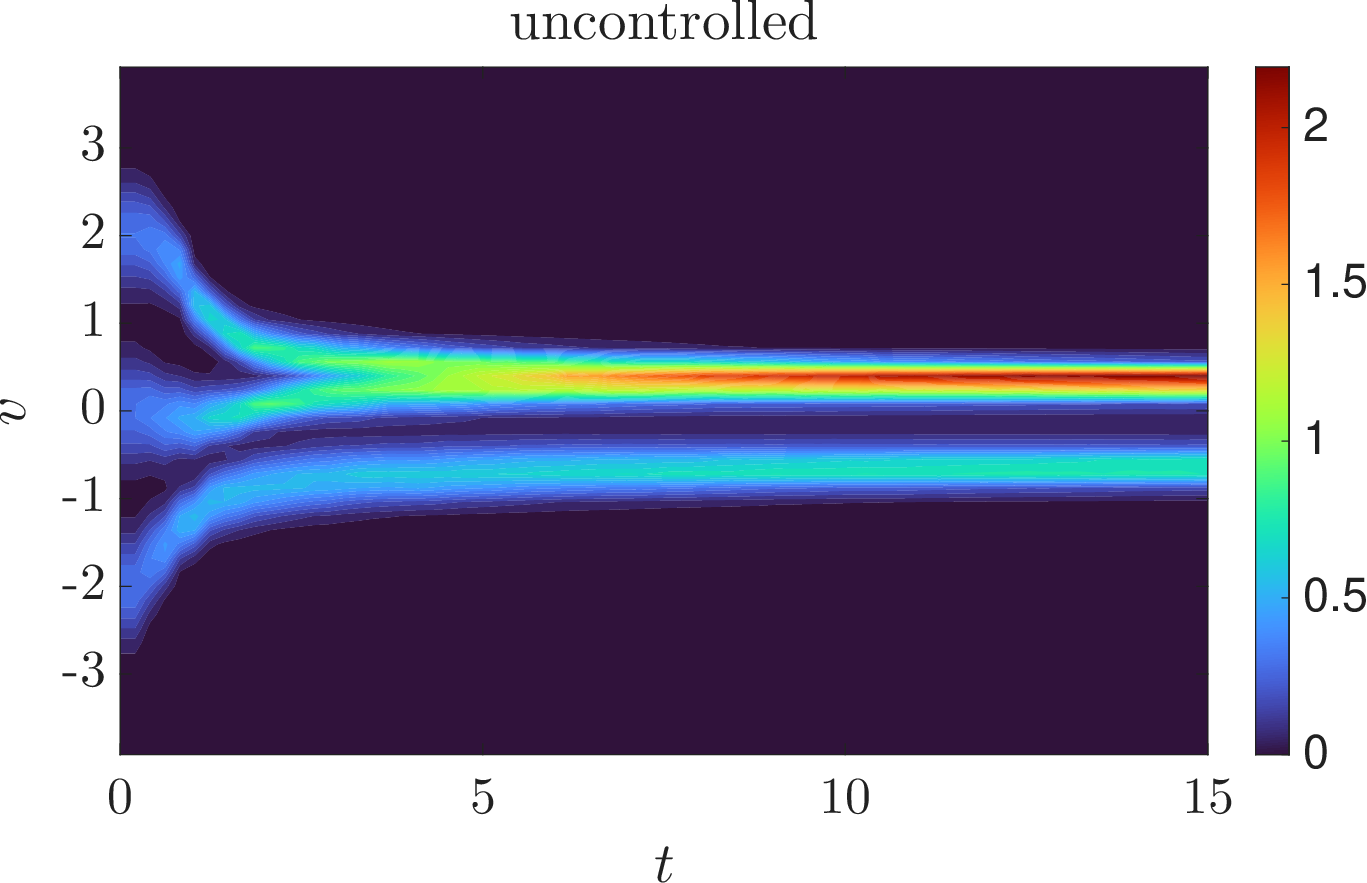}
    \caption{\textbf{Test 2 (uncontrolled).} Mean-field dynamics without control: 
    left panel shows the spatial distribution, and right panel the velocity marginal.}
    \label{fig:meanfield2_uncontrolled}
\end{figure}

\begin{figure}[htbp]
    \centering
    \includegraphics[width=0.485\linewidth]{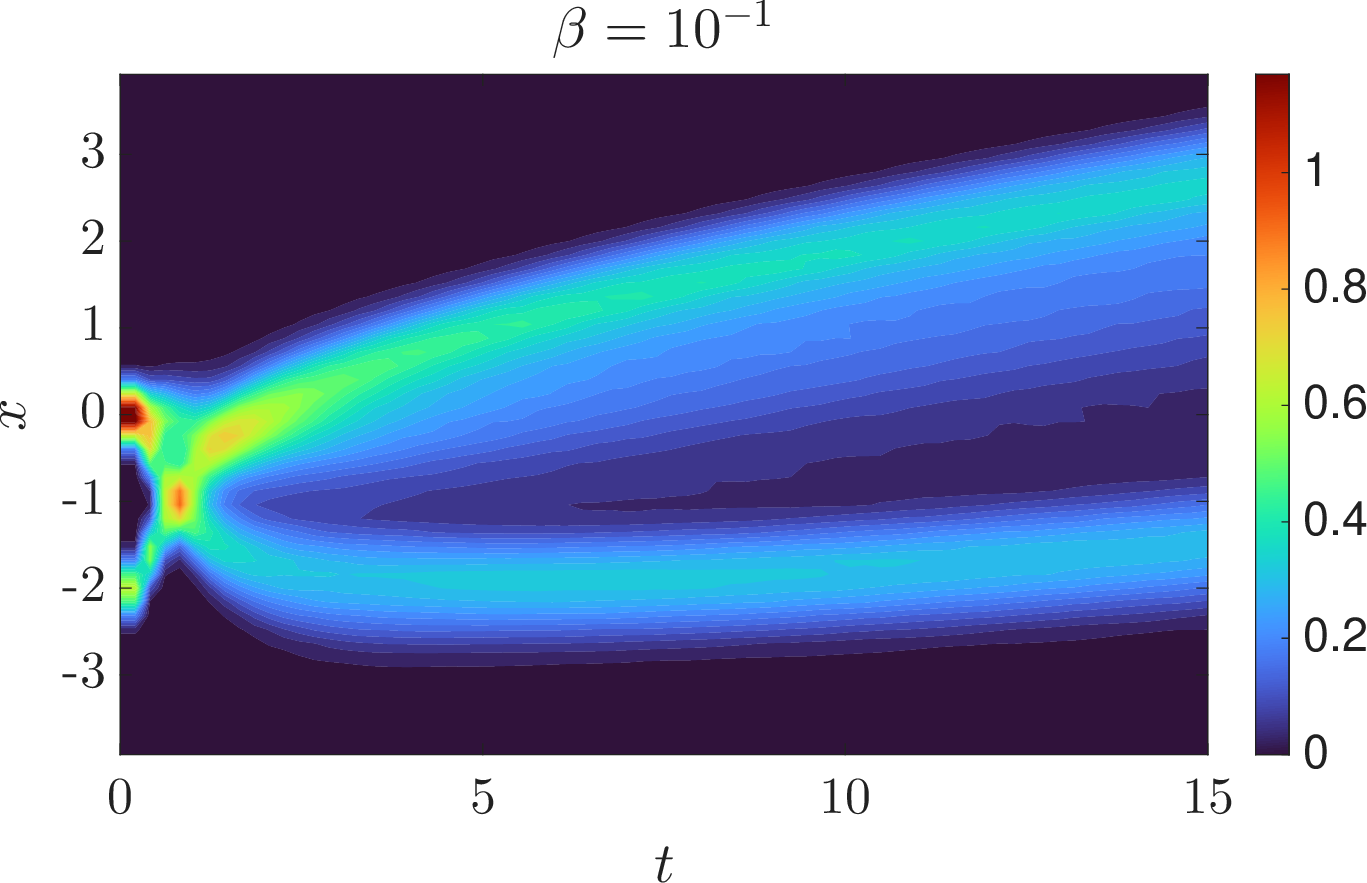}
    \includegraphics[width=0.485\linewidth]{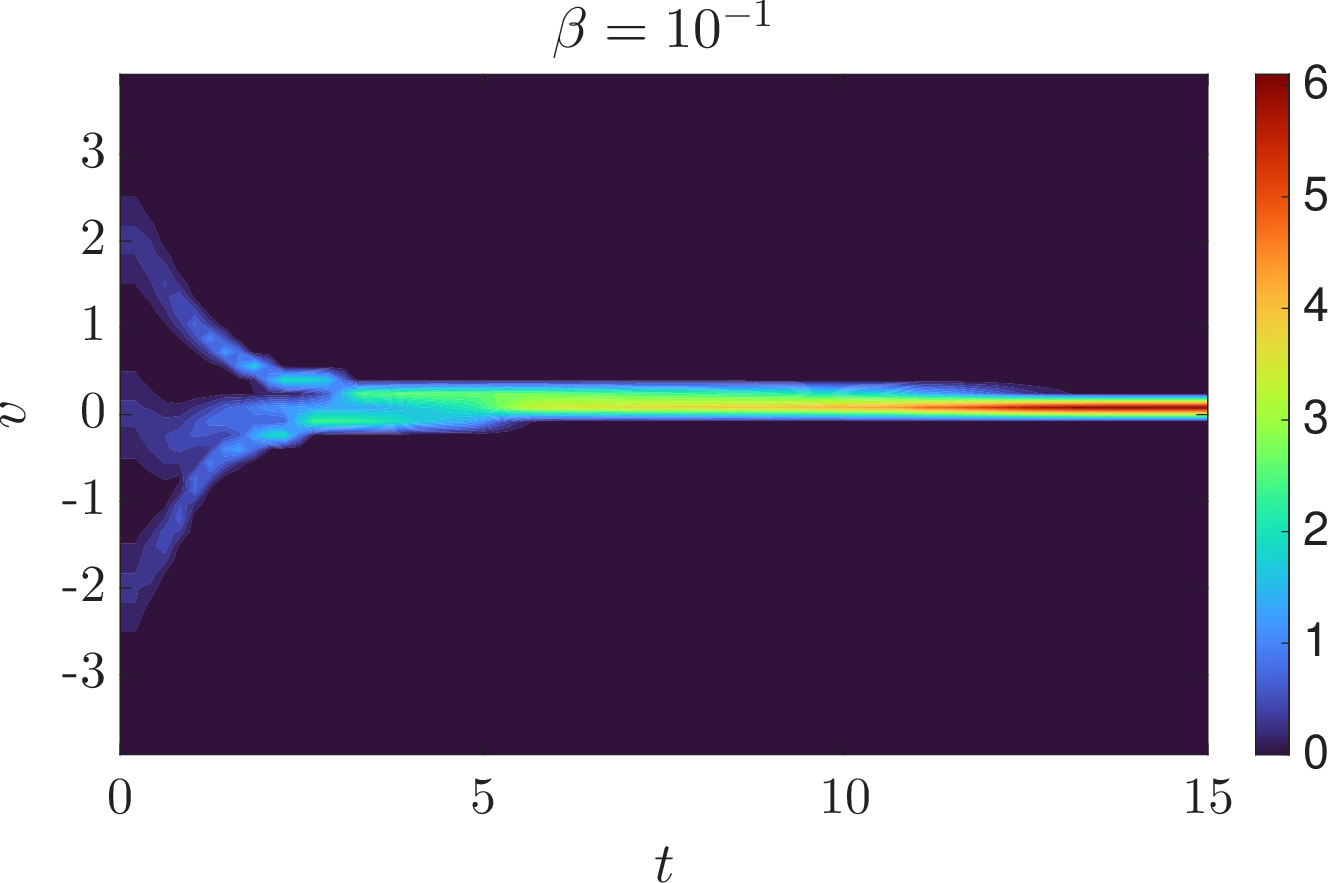}\\[3pt]
    \includegraphics[width=0.485\linewidth]{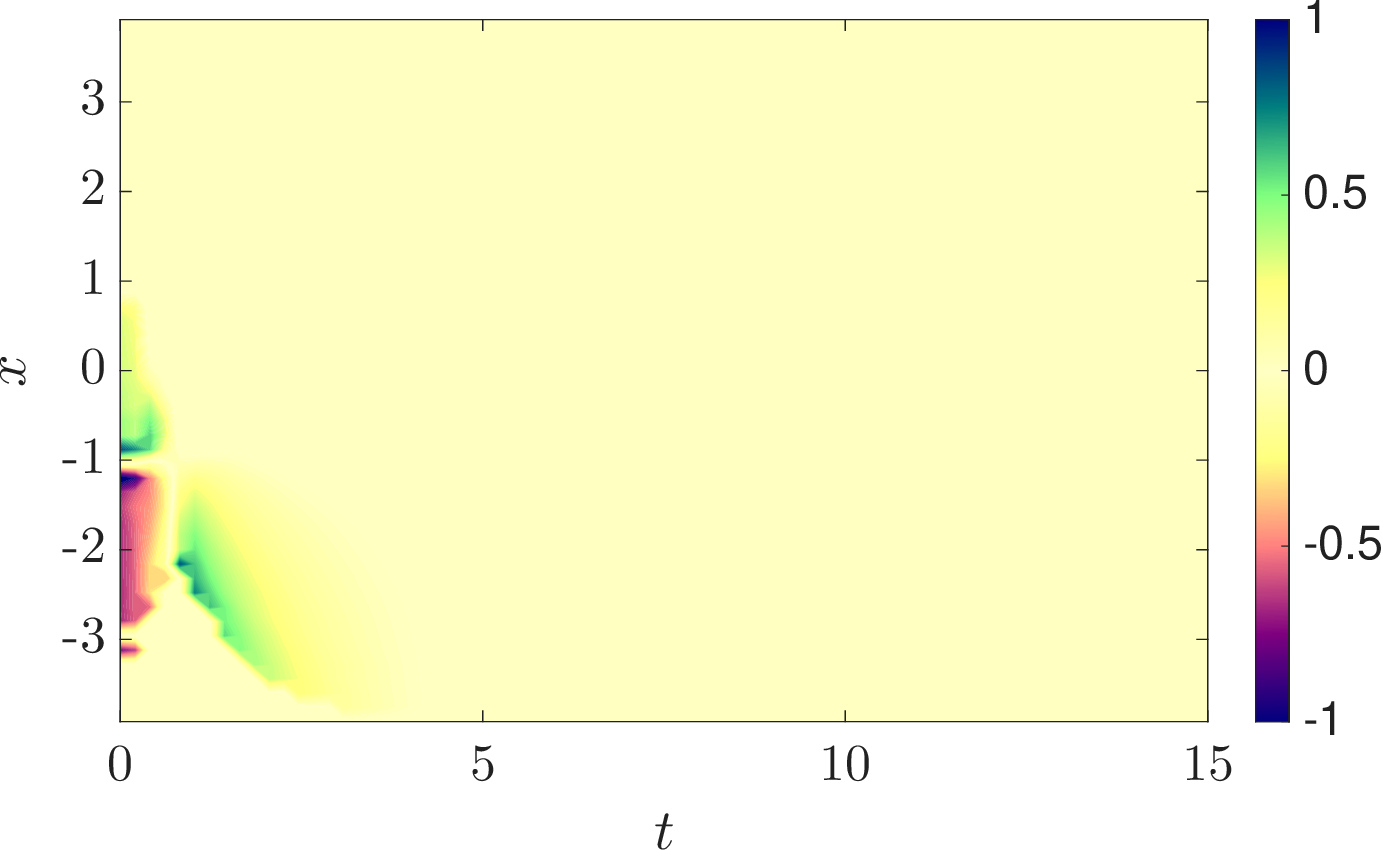}
    \includegraphics[width=0.485\linewidth]{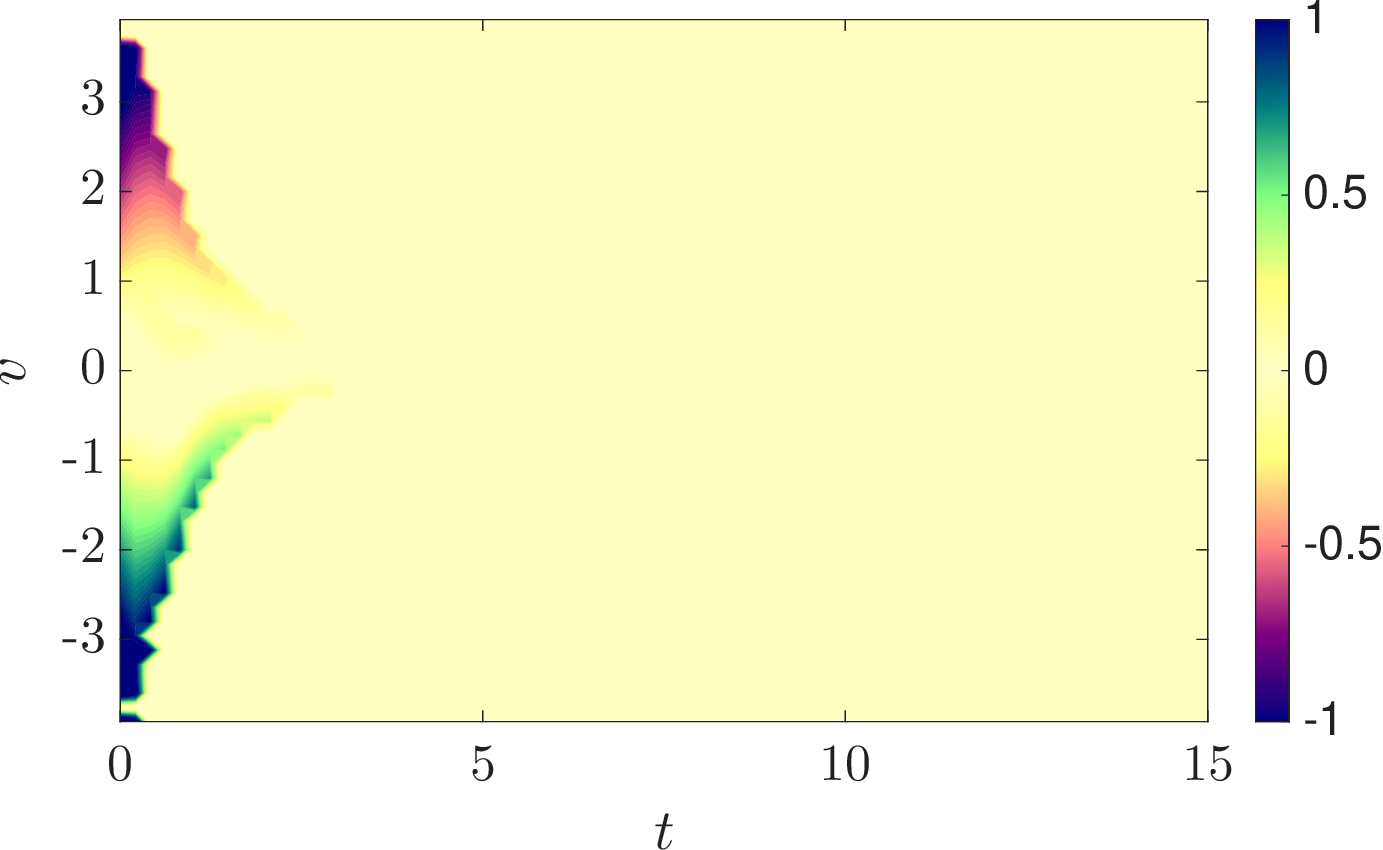}
    \caption{\textbf{Test 2 (controlled).} \emph{Top row:} Mean-field dynamics with control parameters $\alpha =10^{-2}, \ \beta = 10^{-1}$: 
    spatial concentration (left) and velocity alignment (right) are achieved. \emph{Bottom row:} time evolution of control marginals in $x$ (left) and $v$ (right), showing strong activity initially that vanishes as consensus emerges.}
    \label{fig:meanfield2_controlled}
\end{figure}
Figure~\ref{fig:meanfield2_uncontrolled} shows the mean-field evolution in the uncontrolled case, where agents remain dispersed in both space and velocity. Figure~\ref{fig:meanfield2_controlled} reports in the top row the sparse control case where the velocity marginal concentrates, and the spatial distribution is confined. The control efficiently promotes consensus while maintaining sparsity, as depicted in the bottom row for space (left) and velocity (right).

\begin{figure}[htbp]
	\centering
 	\includegraphics[width=0.8\linewidth]{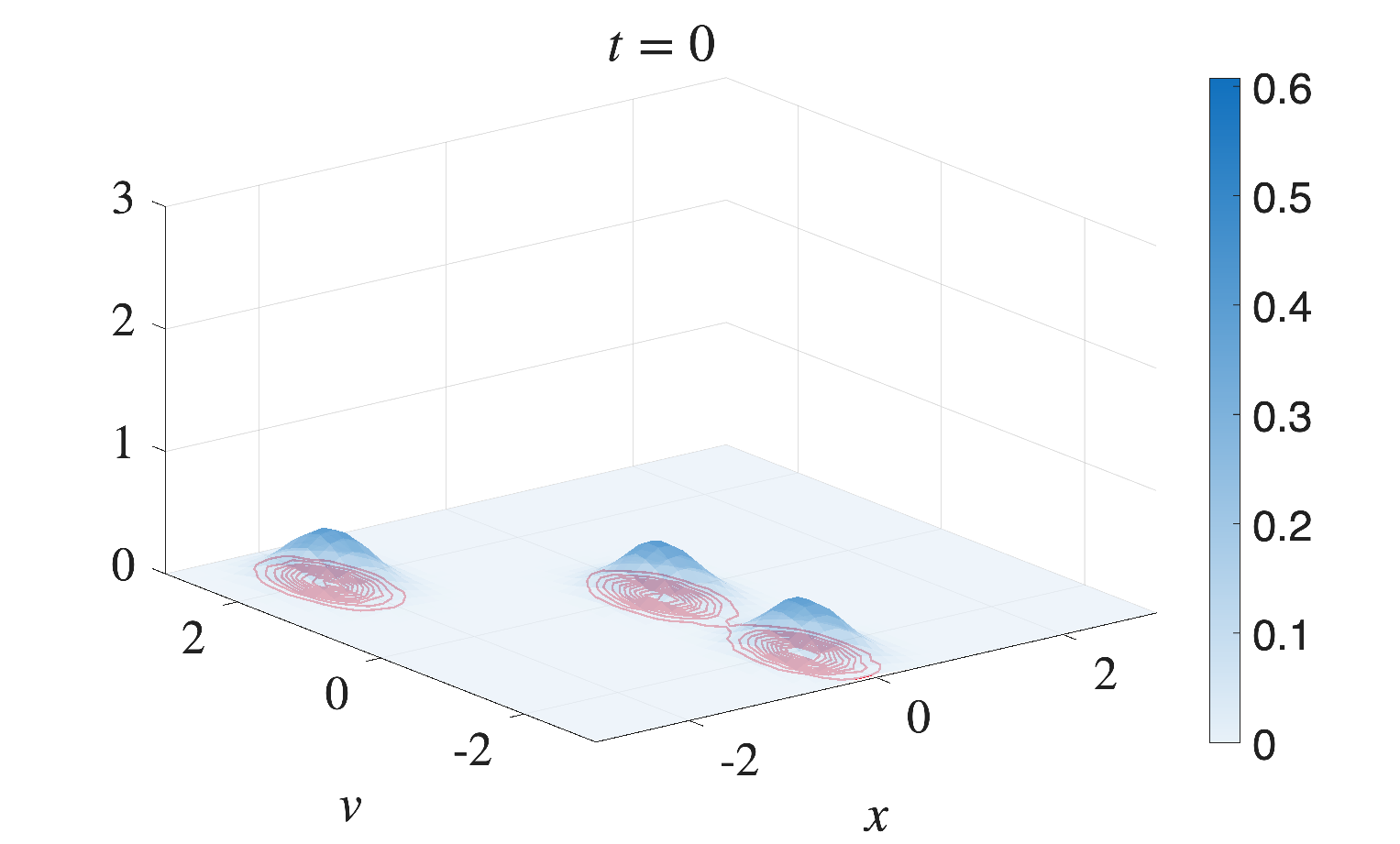}\\
	\includegraphics[width=0.8\linewidth]{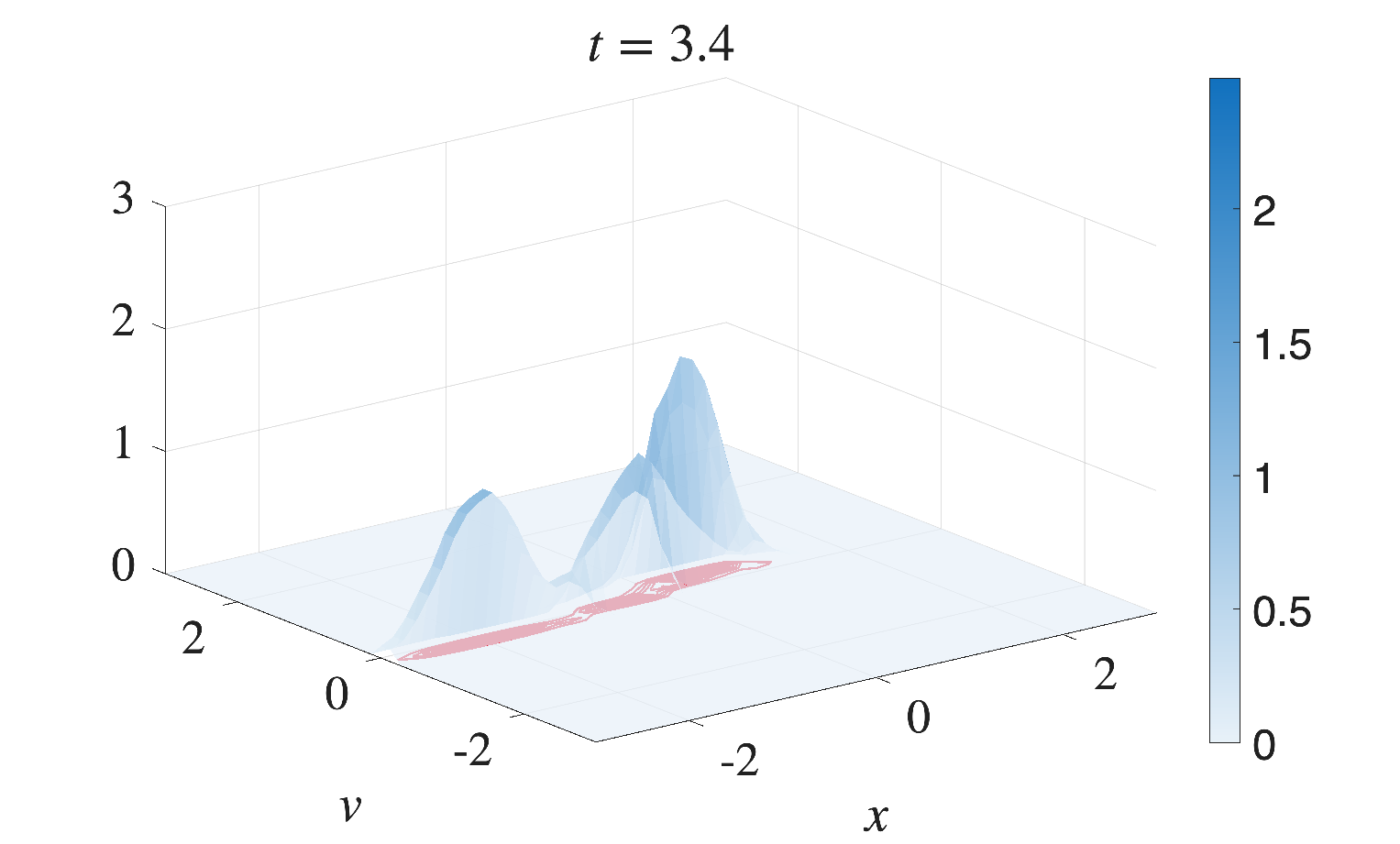}
\caption{\textbf{Test 2.} \emph{Top row:} Initial configuration at $t=0$; \emph{Bottoom row:}  controlled state at $t=3.4$ for $\beta=10^{-1}$.}
	\label{fig:meanfield_snapshots}
\end{figure}

Figure~\ref{fig:meanfield_snapshots} illustrates the phase-space evolution of the mean-field density.
Initially, agents form three clusters that merge into a single aligned group, where here we depict the behavior at time $t=3.4$ where alignment is already almost reached. 
Finally, the control performance is assessed using the Lyapunov functional $V(t)=N_s^{-1}\sum_{i,j}\|v_i-v_j\|^2$, and the number of inactive control components. 
With $N_sN_T=3\times10^6$ control variables, Table~\ref{tab:meanfield_metrics} summarizes the results.

\begin{table}[h!]
\centering
\caption{\textbf{Test 2.} Control performance with $\alpha=0.01$.}
\label{tab:meanfield_metrics}
\begin{tabular}{lcc}
\toprule
Scenario & $V(T)$ & Inactive controls \\
\midrule
Uncontrolled & $2\times10^{-1}$ & $3\times10^6$ \\
Sparse ($\beta=10^{-1}$) & $5\times10^{-4}$ & $2\times10^6$ \\
Full control & $9\times10^{-7}$ & $0$ \\
\bottomrule
\end{tabular}
\end{table}

Sparse control greatly reduces dispersion with two-thirds of the components inactive, while full control achieves near-perfect alignment at full activation, illustrating the performance–sparsity trade-off.

\subsection{Test 2: Mean-field two dimensional dynamics}

\begin{figure}[htbp]
    \centering    
        \includegraphics[width=0.485\linewidth]{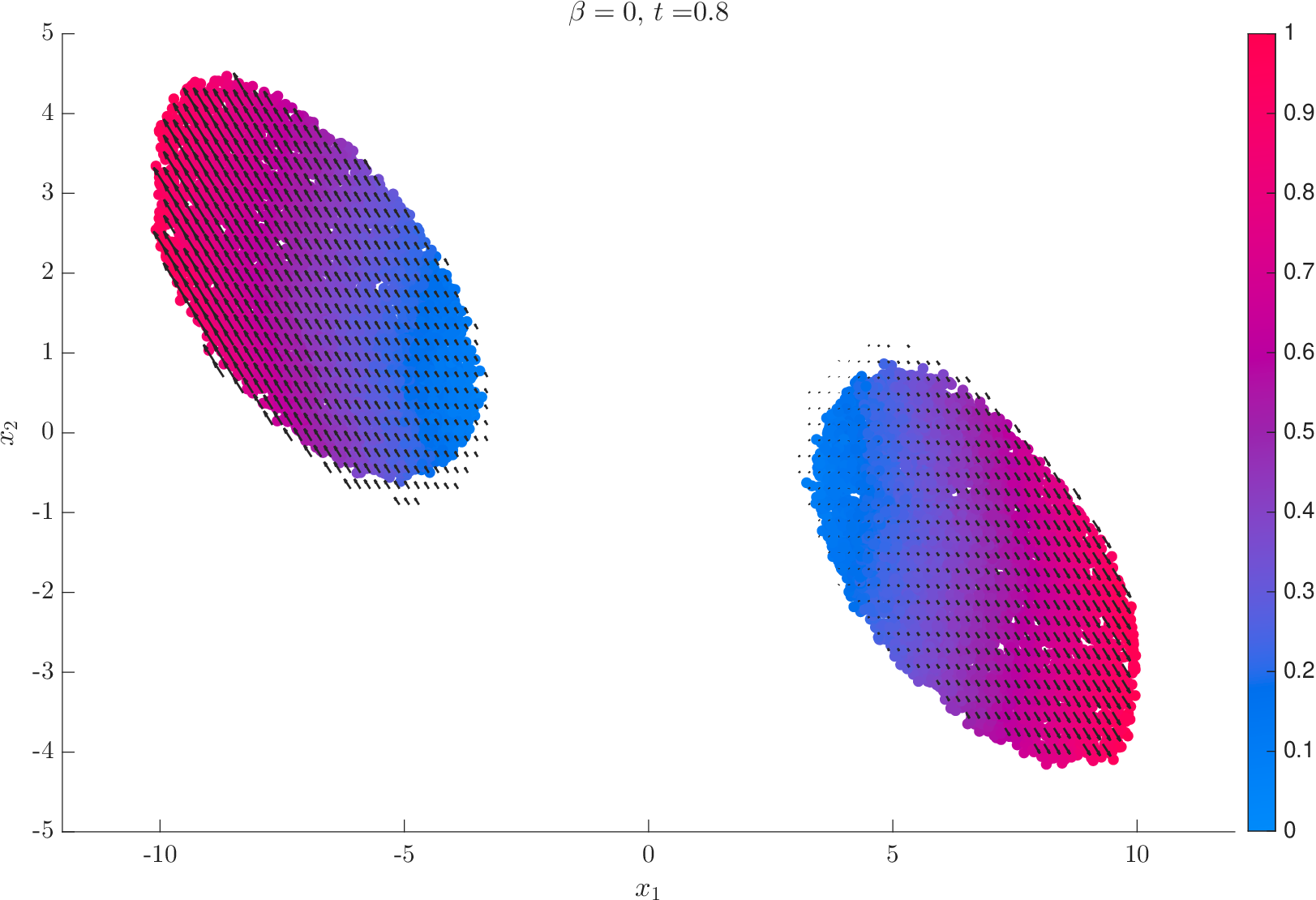}
    \includegraphics[width=0.485\linewidth]{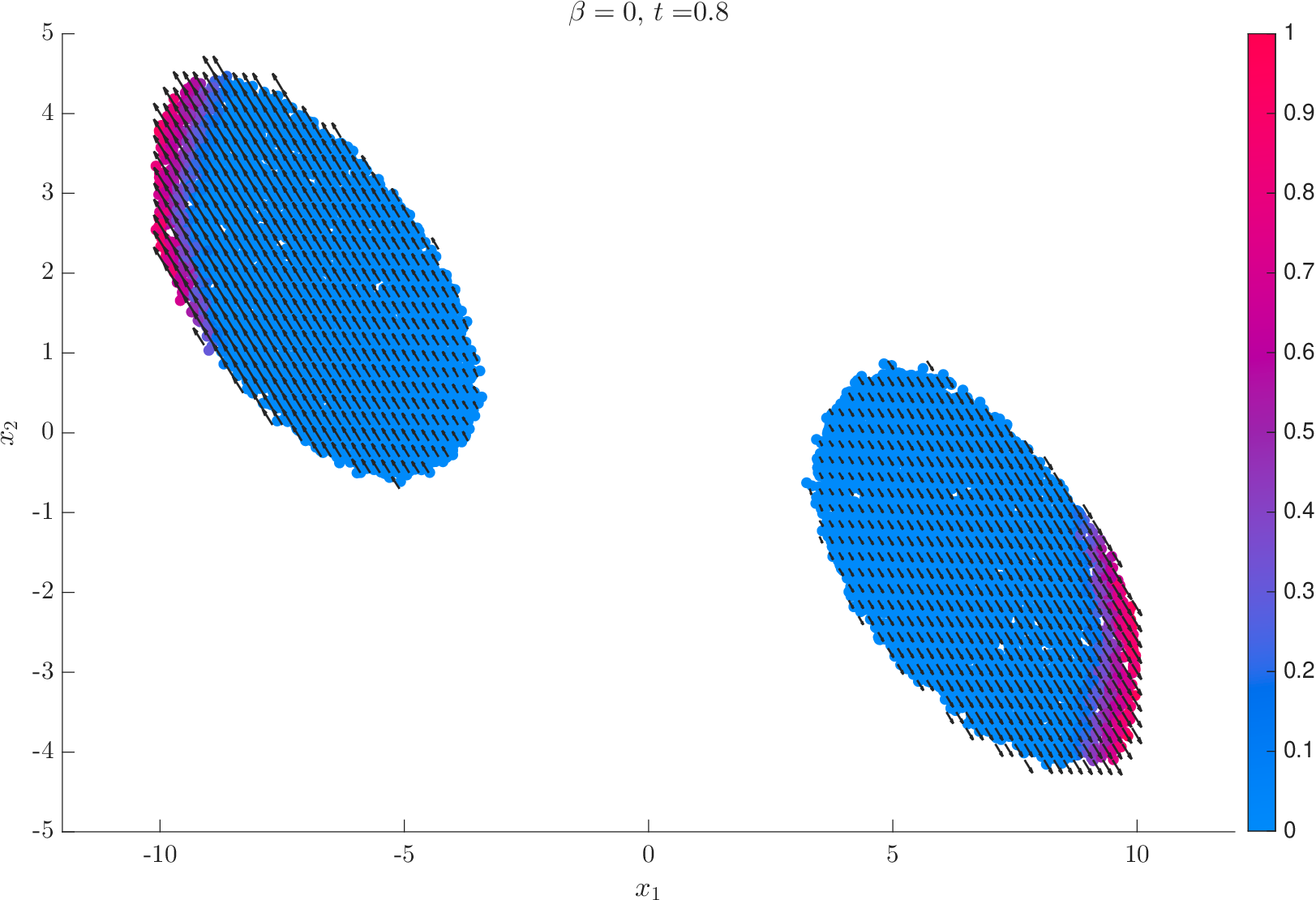}
    \caption{\textbf{Test 3.} 
\emph{Bottom row:} Spatial marginals at $t=0.8$ for $\beta=0$ (left) and $\beta=0.1$ (right). Arrows indicate the momentum field, the scale of color the magnitude of the control.}
    \label{fig:mf_2d_test}
\end{figure}

\begin{figure}[htbp]
    \centering
      \includegraphics[width=0.485\linewidth]{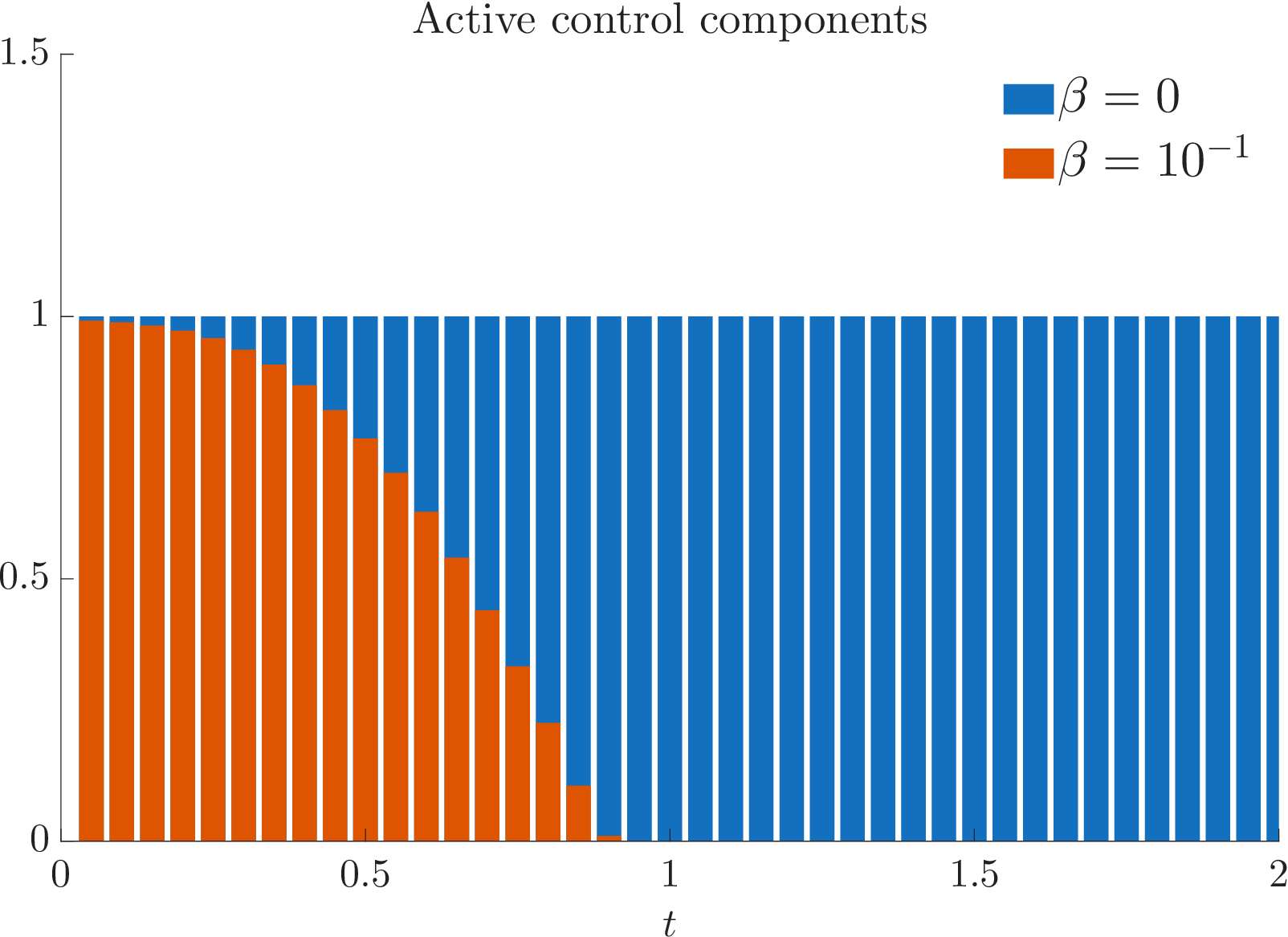}\,
       \includegraphics[width=0.485\linewidth]{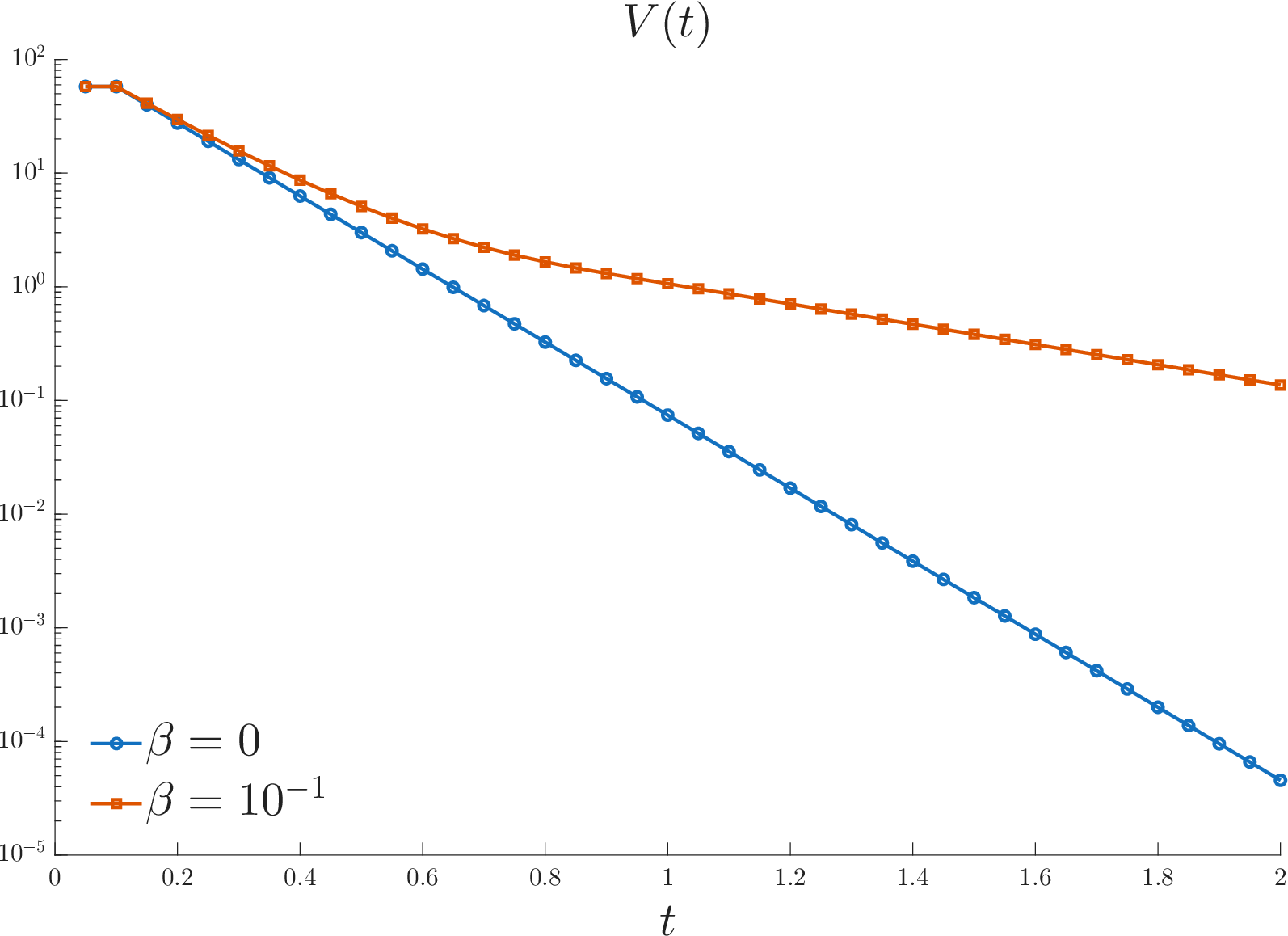}
    \caption{\textbf{Test 2.}  
\emph{Left:} Active control components over time for $\beta=0$ and $\beta=10^{-1}$. \emph{Right:} Decay of $V(t)$ for the different $\beta$ penalization regime.}
    \label{fig:mf_2d_controls}
\end{figure}

We now consider the two-dimensional mean-field CS system with final time $T=2$, time step $\Delta t=0.05$, and $N=10^4$ particles. The initial density $f^0(x,v)$ consists of two symmetric clusters centered at $(-5,0)$ and $(5,0)$, each uniformly distributed over a disc of radius $r_{\max}=2$. Velocities depend quadratically on the distance from the origin, and include small Gaussian noise, producing two groups moving respectively along the directions $(-1,1)$ and $(1,-1)$. This yields a strongly misaligned configuration. The control problem is discretized as in \eqref{eq:disc_particles}–\eqref{eq:J_discrete}, with quadratic penalization parameter $\alpha=0.01$, sparsity $\beta=0.1$, and global budget constraint $\widetilde B=10^6$. The TOS algorithm is run with the same parameters used in the previous test. Figure~\ref{fig:mf_2d_test} shows the two-dimensional mean-field test, reporting the spatial marginals
at $t=0.8$ for $\beta=0$ (left) and $\beta=0.1$ (right).
Without sparsity, the control induces stronger deformation, while for $\beta=0.1$ the action is more localized. Arrows show the momentum field, and the color intensity reflects the magnitude of the control.
~Figure~\ref{fig:mf_2d_controls} shows the activation pattern of the control and the decay of $V(t)$.  
For $\beta=0$ the control stays fully active and $V(t)$ decreases rapidly, whereas $\beta=10^{-1}$ yields sparser activation and a slower decay.

\section{Conclusion}
This work presented an efficient framework for sparse stabilization of large-scale multi-agent systems, balancing control performance and computational cost.
Using proximal methods, we obtained sparse controls that drive the system to consensus with reduced effort, while the mean-field formulation ensured scalability. The numerical results consistently confirm the effectiveness of the approach.

\bibliographystyle{IEEEtran}\bibliography{IEEEabrv,biblioCS2}

\end{document}